\documentclass{aims}
\usepackage{graphicx}  %\usepackage{epstopdf}%This is to transfer .eps figure to .pdf figure; please compile your paper using PDFLeTex or PDFTeXify.

\usepackage{bm}
\usepackage[colorlinks=true]{hyperref}
\usepackage{algorithm}
\usepackage{algorithmic}
\hypersetup{urlcolor=blue, citecolor=red}
\newcommand{\sg}{\mbox{sgn}}

\def\currentvolume{X}
 \def\currentissue{X}
  \def\currentyear{X}
   \def\currentmonth{X}
    \def\ppages{X}
     \def\DOI{X}

\title[Permutation Binary Neural Networks] %Use the shortened version of the full title
{Permutation Binary Neural Networks: Analysis of Periodic Orbits and Its Applications}

\author[Hotaka Udagawa, Taiji Okano, Toshimichi Saito]{}


 \email{tsaito@hosei.ac.jp}

\thanks{$^*$ Corresponding author: Toshimichi Saito}

\begin{document}
\maketitle

\centerline{\scshape Hotaka Udagawa, Taiji Okano and Toshimichi Saito$^*$}
\medskip
{\footnotesize
 \centerline{Department of Electrical and Electronic Engineering, HOSEI University}
   \centerline{Koganei, Tokyo, 184-8584, Japan}
} % Do not forget to end the {\footnotesize by the sign }

\bigskip

%\centerline{(Communicated by the associate editor name)}

\begin{abstract}
This paper presents a permutation binary neural network characterized by 
local binary connection, global permutation connection, and the signum activation function.  
The dynamics is described by a difference equation of binary state variables. 
Depending on the connection, the network generates various periodic orbits of binary vectors. 
The binary/permutation connection brings benefits to precise analysis and to FPGA based hardware implementation. 
In order to consider the periodic orbits, we introduce three tools:  
a composition return map for visualization of the dynamics,  
two feature quantities for classification of periodic orbits, and 
an FPGA based hardware prototype for engineering applications. 
Using the tools, we have analyzed all the 6-dimensional networks. 
Typical periodic orbits are confirmed experimentally.
%
%\textbf{200} words.
\end{abstract}

\section{Introduction}
\label{intro}
Discrete-time recurrent-type neural networks have been studied actively 
in both viewpoints of nonlinear dynamics and engineering applications \cite{hpf}-\cite{adachi}. 
The networks are characterized by real valued connection parameters and nonlinear activation functions. 
The dynamics is described by autonomous difference equations of real valued state variables.     
Even in the low-dimensional cases, 
the networks can exhibit various periodic/chaotic phenomena \cite{adachi}-\cite{ott}. 
Several effective synthesis methods have been presented for applications such as 
associative memories \cite{hpf}-\cite{adachi}, 
combinatorial optimization problems solvers \cite{tsp}, and 
reservoir computing \cite{rc1}. 
In the recurrent-type neural networks, fixed points have been studied sufficiently. 
For example, parameter conditions for existence/stability of fixed points are presented 
and are applied to storage of stable desired memories into the associative memories \cite{hpf}-\cite{am2}. 
However, analysis of periodic orbits is difficult and has not been studied sufficiently. 
The difficulty is caused by 
huge parameter space, curse of dimensionality and complex nonlinear dynamics. 
In order to realize the analysis, 
a simple mathematical model with a variety of periodic orbits is necessary. 

This paper presents a permutation binary neural network (PBNN) 
and considers classification/stability of periodic orbits. 
The PBNN is a simple recurrent-type neural network characterized by 
local binary connection parameters, global permutation connection parameters, and the signum activation function.  
The dynamics is described by an $N$-dimensional autonomous difference equation of binary state variables. 
Depending on the parameters and initial conditions, the PBNN generates various periodic orbits of binary vectors (BPOs). 
Engineering applications of the BPOs include 
control signal of switching power converters \cite{pe} and  
control signal of walking robots \cite{cpg1}-\cite{cpg3}. 
The PBNNs can be regarded as a simplified version of dynamic binary neural networks (DBNNs \cite{dbnn1}-\cite{cnna}) characterized by global ternary connection. 
If the permutation connection does not exist, 
the PBNNs are equivalent to elementary cellular automata \cite{ca1}-\cite{error} governed by rules of linearly separable Boolean functions. 
The binary/permutation connection brings benefits to precise analysis of dynamics and 
to simple hardware implementation based on a field-programmable gate array (FPGA) \cite{fpga} \cite{dsn}. 

In order to consider the PBNNs and their BPOs, we present three tools. 
The first tool is a composition return map (Cmap) for visualization of the dynamics. 
The Cmap is a composition of two maps. 
The first map represents the local binary connection and the second map represents the global permutation connection. 
The second tool is two feature quantities for classification of PBNNs. 
Roughly speaking, the first feature quantity evaluates complexity of the BPOs 
and the second feature quantity evaluates stability of the BPOs. 
Using the two feature quantities, a feature plane is constructed on which PBNNs are classified based on their BPOs. 
The third tool is an FPGA based hardware prototype for engineering applications. 
Using the tools, we have considered all the 6-dimensional PBNNs 
whose parameter space consists of ($8 \times 6!$) points.  
The dynamics is visualized in Cmap on $2^6$ points where we can grasp BPOs and transient phenomena. 
The BPOs from the PBNNs are classified on the feature plane 
where we can see that the global permutation connection makes various BPOs. 
Using the standard design procedure, 
the PBNNs are implemented on an FPGA board and typical BPOs are confirmed experimentally. 
These results provide basic information to realize more detailed analysis of PBNNs and its applications. 

It should be noted that this is the first paper of the PBNN and analysis of its BPOs.  
The Cmap is a developed version of the digital return maps \cite{dbnn1} for DBNNs without the concept of composition. 
Design of the FPGA based hardware prototype is based on simple Boolean functions and is applicable to hardware implementation of the elementary cellular automata.

\section{Two related systems}
As background of the PBNNs, we introduce 
elementary cellular automata (ECAs \cite{ca1}-\cite{rc2}) and 
dynamics binary neural networks (DBNNs, \cite{dbnn1}-\cite{cnna}). 
In various discrete-time dynamical systems, these two systems are simple and exhibit various BPOs. 

\subsection{Elementary cellular automata}
\begin{figure}[tb]
\centering
\includegraphics[width=0.6\columnwidth]{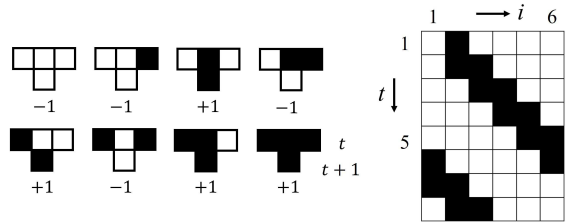}
\caption{
An example of 6-dimensional ECA (RN212) and a sequence of binary vectors. 
Black and white squares denote output $+1$ and $-1$, respectively.  
}
\label{fg1} 
\end{figure}

We introduce the ECAs defined on a ring of $N$ cells. 
Let $x_i^t \in \{ -1, +1 \} \equiv \bm{B}$
\footnote{
this expression is used instead of expression 
$x_i^t \in \{0, 1\}$ in \cite{ca1} \cite{ca2} for matching with PBNNs.
}
be the $i$-th binary state at discrete time $t$. 
The dynamics is described by the $N$-dimensional autonomous difference equation:
\begin{equation}
x^{t+1}_i = F(x^t_{i-1},x^t_i,x^t_{i+1}), \ \ i \in \{ 1, \cdots, N \}, \ N \ge 3
\end{equation}
where $x_0^t \equiv x_N^t$ and $x_{N+1}^t \equiv x_1^t$ for ring-type connection.  
A Boolean function $F$ transforms three binary inputs into one binary output. 
Fig. \ref{fg1} shows an example:  
\[
\begin{array}{llll}
F(-1, -1, -1) = -1 & F(-1, -1, 1) = -1 & F(-1, 1, -1) = 1 &  F(-1, 1, 1) = -1\\
F(1, -1, -1) =  1 & F( 1, -1, 1) = -1 & F( 1, 1, -1) = 1 &  F(1, 1, 1) = 1 
\end{array}    
\]
Decimal expression of the 8 outputs is referred to as the rule number (RN). 
In this example, the 8 outputs $(-1-1+1-1+1-1+1+1) \equiv (00101011)_2 = 212_{10}$ is declared as RN212. 
There exist $2^{2^3} = 256$ rules in the ECAs and 
the rules are characterized by the $\lambda-$parameter = (the number of output +1)/8. 
The RN212 is characterized by $\lambda = 0.5$. %that can give complicated binary sequences.  
Depending on the rules and initial conditions, the ECAs can 
generate various sequences of binary vectors (BPOs and transient phenomena) as studied in \cite{ca1}-\cite{ca3}. 
Engineering applications are many, including 
information compression \cite{ca4}, reservoir computing \cite{rc2}, and error correcting code \cite{error}.   
Analysis of the ECAs is important from both viewpoints of nonlinear dynamics and engineering applications.   

\subsection{Dynamic binary neural networks}
The DBNNs are recurrent-type neural networks 
described by the $N$-dimensional autonomous difference equation of binary state variables:
\begin{equation}
\begin{array}{c}
\displaystyle x_i^{t+1} = \sg \left(\sum_{j=1}^N w_{ij} x_j^t - T_i \right) 
\ i \in \{ 1, \cdots, N \}, \ N \ge 3\\
\sg(X) = \left\{\begin{array}{ll}
+1 & \mbox{for } X \ge 0\\
-1 & \mbox{for } X < 0
\end{array}\right.
\end{array}
\label{dbnn}
\end{equation}
where 
$x_i^t \in \bm{B}$ is the $i$-th binary state at discrete time $t$. 
The DBNNs are characterized by ternary connection parameters $w_{ij} \in \{-1, 0, +1\}$ 
and the signum activation function. 
The threshold parameters $T_j$ are integer. 
Depending on the parameters and initial conditions, the DBNNs can generate various BPOs. 
In \cite{dbnn1}, we have considered storage and stabilization of special BPOs 
for application to switching power converters. 
In \cite{dbnn2}, we have analyzed the DBNNs with special sparse connection and 
have given theoretical results for storage and stabilization of rotation-type BPOs. 
However, in these DBNNs, storage of any BPO is not guaranteed. 
In \cite{dbnn3}, we have presented three-layer DBNNs and 
a simple parameter setting method that guarantees storage of any BPO. 
However, the number of hidden neurons increases in proportional to the period of BPOs. 
The three-layer DBNNs have been implemented in an FPGA board and have been applied to 
controlling hexapod walking robots \cite{dbnn4} \cite{cnna}. 

\subsection{Problems in ECAs and DBNNs}
Although BPOs are simpler than real valued periodic orbits, it is not easy to analyze BPOs in DBNNs and ECAs. 
Connection parameter space of the DBNNs is huge. 
The ECAs and DBNNs exhibit complicated dynamics as $N$ increases.

\section{Permutation Binary Neural Networks and Identifier} 

\begin{figure}[tb]
\centering
\includegraphics[width=0.7\columnwidth]{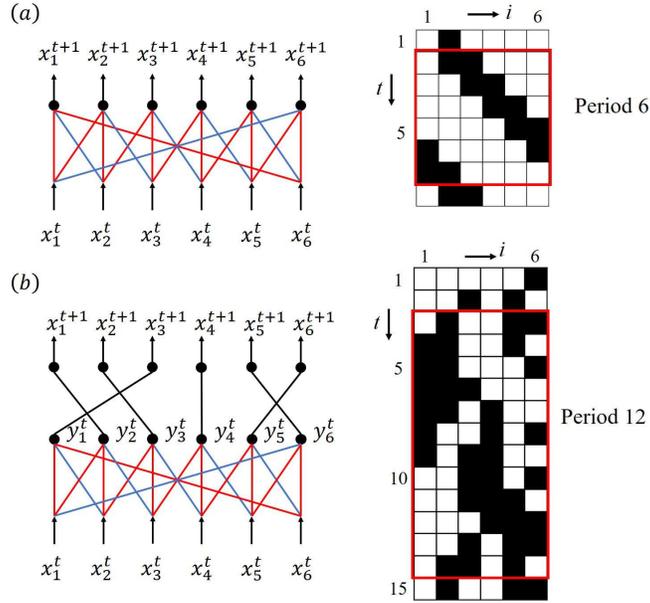}
\caption{
Networks and BPOs. 
(a) SBNN in Example 1. 
Red and blue branches denote 
positive and negative connections, respectively. 
(b) PBNN in Example 2 (CN6, P231465).  
Black branches denote permutation connection. 
}
\label{fg2} 
\end{figure}

We presents the PBNN and consider BPOs from the PBNN. 
First, we introduce a simple binary neural network (SBNN)  
corresponding to the input to hidden layers in the PBNN.  
The SBNN is characterized by local binary connection and the signum activation function. 
The dynamics is described by the following 
$N$-dimensional autonomous difference equation of binary state variables:
\begin{equation}
x_i^{t+1} = \sg \left(w_a x_{i-1}^{t} + w_b x_i^t + w_c x_{i+1}^t \right), \ i \in \{1, \cdots, N \}, \ N \ge 3
\end{equation}
where $x_i^t \in \bm{B}$ is the $i$-th binary state variable at discrete time $t$. 
As shown in Fig. \ref{fg2}(a), 
$x_{N+1}^t \equiv x_1^t$ and $x_0^t \equiv x_N^t$ for ring-type connection.  
The local binary connection is determined by 
three binary valued connection parameters $w_a \in \bm{B}$, $w_b \in \bm{B}$, and $w_c \in \bm{B}$. 
Depending on the parameters, the signum activation function $\sg(X)$ constructs a linearly separable Boolean function  
from three binary inputs $(x^t_{i-1}, x^t_i, x^t_{i+1})$ to one binary output $x_i^{t+1}$. 
Note that $X \in \{-3, -1, +1, +3\}$ for $X \equiv w_a x_{i-1}^{t} + w_b x_i^t + w_c x_{i+1}^t$. 
In order to identify the local binary connection, we introduce connection number (CN) 
\begin{equation}
\begin{array}{cc}
\mbox{CN0 for } (w_a, w_b, w_c) = (-1, -1, -1) & \mbox{(RN23)} \\
\mbox{CN1 for } (w_a, w_b, w_c) = (-1, -1, +1) & \mbox{(RN43)} \\
\mbox{CN2 for } (w_a, w_b, w_c) = (-1, +1, -1) & \mbox{(RN77)} \\
\mbox{CN3 for } (w_a, w_b, w_c) = (-1, +1, +1) & \mbox{(RN142)} \\
\mbox{CN4 for } (w_a, w_b, w_c) = (+1, -1, -1) & \mbox{(RN113)} \\
\mbox{CN5 for } (w_a, w_b, w_c) = (+1, -1, +1) & \mbox{(RN178)} \\
\mbox{CN6 for } (w_a, w_b, w_c) = (+1, +1, -1) & \mbox{(RN212)} \\
\mbox{CN7 for } (w_a, w_b, w_c) = (+1, +1, +1) & \mbox{(RN232)}
\end{array}
\label{cn}
\end{equation}
The 8 CNs corresponding to 8 Boolean function of 8 RNs of ECAs as shown in Eq. (\ref{cn}). 
Fig. \ref{fg2} (a) shows Example 1: A 6-dimensional SBNN of CN6. 
Depending on initial conditions, this SBNN exhibits sequences of binary vectors 
including a BPO with period 6 as shown in the figure. 
All the 8 rules are characterized by $\lambda$-parameter=0.5. 

%%%%%%%%%%%%%%%%%%%

%

Adding global permutation connection from hidden to output layer, the PBNN is constructed. 
The dynamics is described by the following 
$N$-dimensional autonomous difference equation of binary state variables:
\begin{equation}
\begin{array}{l}
x_i^{t+1} = y_{\sigma(i)}^t\\
y_i^t = \sg \left(w_a x_{i-1}^t + w_b x_i^t + w_c x_{i+1}^t \right)
\end{array} \ \ 
\sigma = \left(
    \begin{array}{cccc}
      1         & 2         & \cdots & N \\
      \sigma(1) & \sigma(2) & \cdots & \sigma(N)
    \end{array}
  \right)
%\end{array}
\label{PBNN}
\end{equation}
where 
$y_i^t \in \bm{B}$ is the $i$-th binary hidden state and $\sigma$ is a permutation. 
Let $\bm{x}^t \equiv (x_1^t, \cdots, x_N^t)$ and let $\bm{y}^t \equiv (y_1^t, \cdots, y_N^t)$. 
The local binary connection (from input to hidden layers) transforms an binary input vector $\bm{x}^t$ into the binary hidden vector $\bm{y}^t$. 
The global permutation connection (from hidden to output layers) transforms the $\bm{y}^t$ into the binary output vector $\bm{x}^{t+1}$. 
Fig. \ref{fg2} (b) shows Example 2: A 6-dimensional PBNN given by
\begin{equation}
%\begin{array}{l}
(w_a, w_b, w_c) = (+1, +1, -1) \ \ \ \mbox{(CN6)}, \ 
\sigma = \left(
    \begin{array}{cccccc}
      1 & 2 & 3 & 4 & 5 & 6\\
      2 & 3 & 1 & 4 & 6 & 5
    \end{array}
  \right)   \ \ \mbox{(P231465)}
%\end{array}
\end{equation}
As shown in the figure, this PBNN exhibits a BPO with period 12 where
a BPO with period $p$ is defined by 
\begin{equation}
\bm{z}^1, \cdots, \bm{z}^p, \cdots \ \ 
\left\{
\begin{array}{ll}
\bm{z}^{t_1} =   \bm{z}^{t_2} & \mbox{for } |t_2 - t_1| = np\\ 
\bm{z}^{t_1} \ne \bm{z}^{t_2} & \mbox{for } |t_2 - t_1| \ne np 
\end{array}\right. 
\label{bpo}
\end{equation}
where $\bm{z}^t = (z_1^t, \cdots, z_N^t)$, $z_i^t \in \bm{B}$ and $n$ denotes positive integers. 
In Example 2, the local binary connection is the same as the connection of SBNN in Example 1: 
applying the global permutation connection to the SBNN of BPO with period 6, 
we obtain the PBNN that generates BPO with period 12. 
Effects of the permutation connection is discussed in Section \ref{nanalysis}. 
In order to identify the global permutation connection, 
we introduce permutation identifier P$\sigma(1) \cdots \sigma(N)$. 
In Example 2, the PBNN is identified by CN6 and P231465. 
In general, a PBNN is identified by connection number the CN and the permutation identifier. 
Here, we have described the parameter space of a PBNN clearly: 
the parameter subspace of a SBNN consists of 8 points represented by the 8 connection numbers (CNs), 
the parameter subspace of the permutation consists of $N!$ permutation identifiers, and 
the parameter space of a PBNN consists of the $8 \times N!$ points 
represented by the CN and the permutation identifier.
The binary/permutation connection parameters 
bring benefits in precise analysis %by integer arithmetic 
and simple FPGA based hardware implementation.

\section{Tools for numerical analysis}
\label{tools}

\subsection{Tool 1: Composition digital return map}
In order to visualize the dynamics of PBNNs, we present the Tool 1: composition digital return map (Cmap). 
The Cmap is a composition of the 1st map $f_1$ for local binary connection (SBNN) and 
the 2nd map $f_2$ for global permutation connection. 
First, we introduce the 1st map $f_1$ for the SBNN. 
The domain and range of the SBNN are the set of all the $N$-dimensional binary vectors $\bm{B}^N$. 
The $\bm{B}^N$ is equivalent to a set of $2^N$ points $L_N \equiv \{ c_1, \cdots, c_{2^N} \}$, 
e. g., $C_1 \equiv (-1, \cdots, -1)$, $C_2 \equiv (+1. -1. \cdots, -1)$, $\cdots$, $C_{2^N} \equiv (+1, \cdots, +1)$. 
The dynamics of the SBNN is integrated into 
\begin{equation}
\mbox{1st map: }   \bm{x}^{t+1} = f_1(\bm{x}^t), \  \bm{x}^t \in \bm{B}^N \equiv L_D
\end{equation}
Here we give basic definition for the 1st map. 

{\bf Definition 1}: 
A point $\bm{z}_p \in L_D$ is said to be a binary periodic point (BPP) with period $p$ 
if $f^p(\bm{z}_p) = \bm{z}_p$ and $f_1(\bm{z}_p)$ to $f_1^p(\bm{z}_p)$ are all different  
where $f_1^k$ is the $k$-fold composition of $f_1$.  
A sequence of the BPPs, $\{ f_1(\theta_p), \cdots, f_1^p(\theta_p) \}$, 
is said to be a BPO with period $p$. 
A point $\bm{z}_e$ is said to be an eventually periodic point (EPP)
if $\bm{z}_e$ is not a BPP but falls into a BPO, i.e., 
there exists some positive integer $l$ such that $f_1^l(\bm{z}_e)$ is a BPP. 

\begin{figure}[tb]
\centering
\includegraphics[width=0.8\columnwidth]{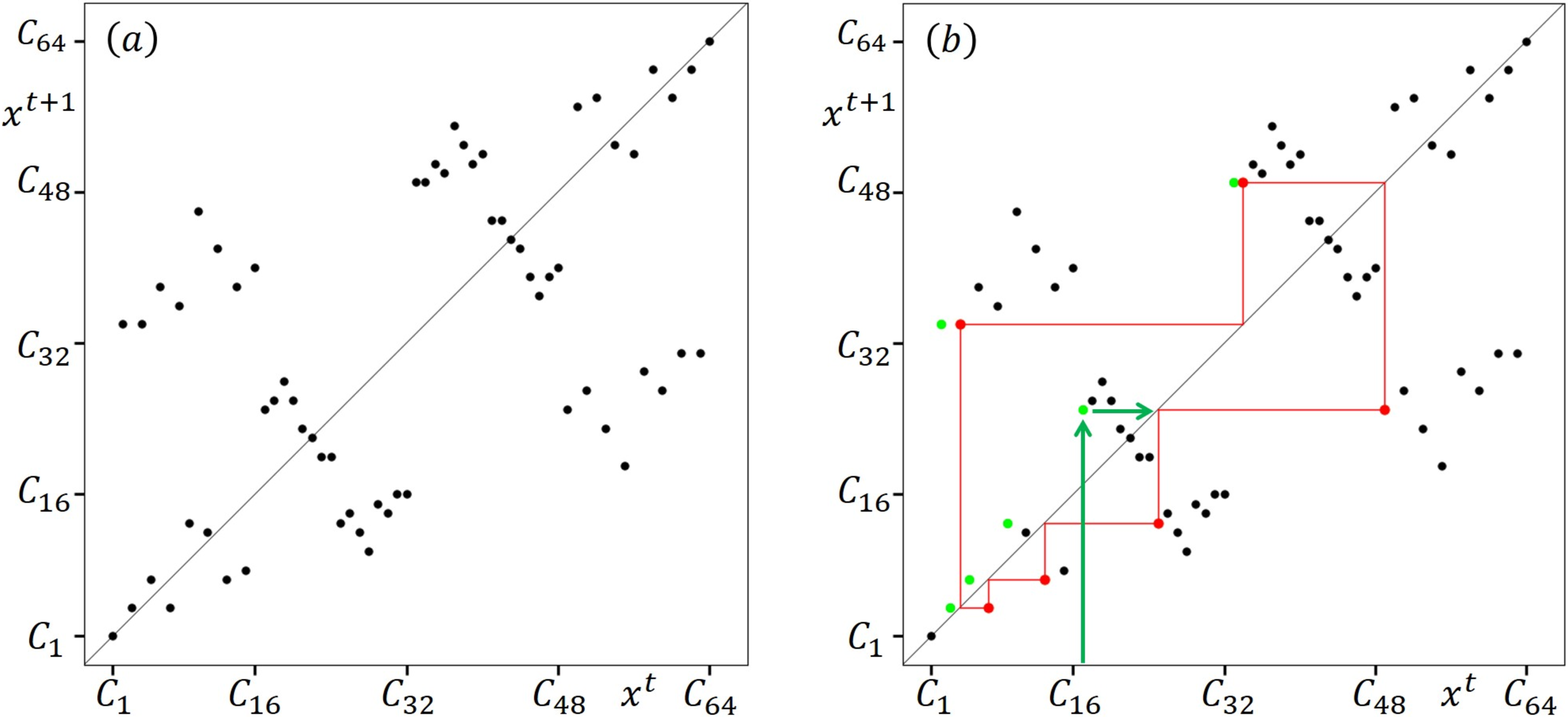}
\caption{1st map $f_1$ for SBNN, CN6. 
(a) 1st map of 64 points. 
(b) BPO(red) with period 6 and EPPs (green).
}
\label{fg3} 
\end{figure}

Fig. \ref{fg3} (a) shows the 1st map consisting of $2^6$ points for the 6-dimensional SBNN in Example 1. 
The 1st map exhibits a BPO with period 6 as shown in Fig. \ref{fg3} (b) 
corresponding to the BPO with period 6 in Fig. \ref{fg2} (a). 
In the figure, several EPPs are also shown.   

\begin{figure}[tb]
\centering
\includegraphics[width=0.8\columnwidth]{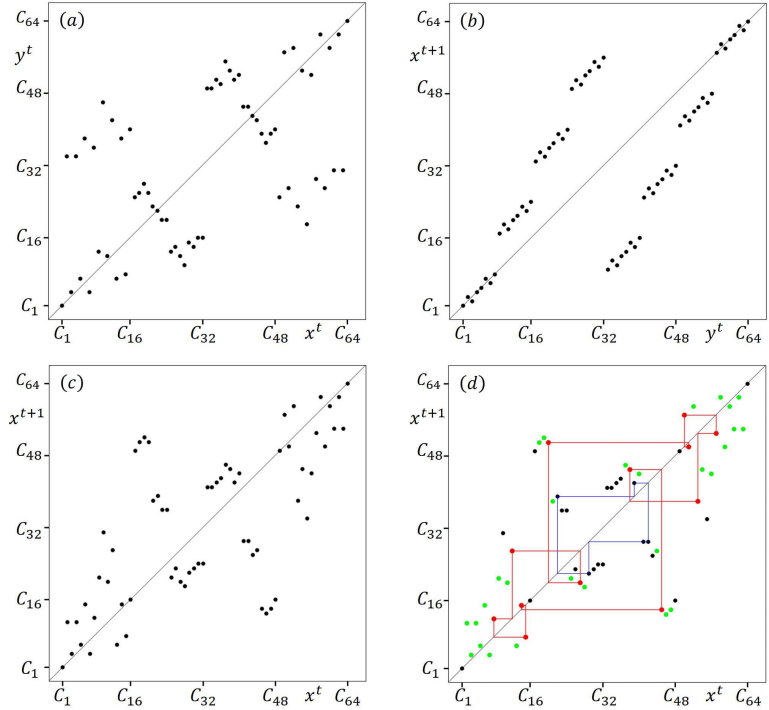}
\caption{Cmaps for PBNN.  
(a) The 1st map $f_1$ for CN6. 
(b) The 2nd map $f_2$ for P231465. 
(c) Cmap $f$ of 64 points. 
(d) BPO with period 12 (red) and EPPs (green). 
BPO with period 4 (blue). 
}
\label{fg4} 
\end{figure}

Next, we introduce the 2nd map $f_2$. 
In the PBNN, the binary hidden state vector $\bm{y}^t \in \bm{B}^N$ 
is transformed into the binary output vector $\bm{x}^{t+1} \in \bm{B}^N$ via permutation $\sigma$. 
The transformation is integrated into 
\begin{equation}
\mbox{2nd map: }   \bm{x}^{t+1} = f_2(\bm{y}^t), \  \bm{y}^t \in \bm{B}^N \equiv L_D. 
\end{equation}
The dynamics of the PBNN is integrated into a composition map (Cmap) of the 1st and the 2nd maps: 
\begin{equation}
\mbox{Cmap: }   \bm{x}^{t+1} = f(\bm{x}^t) \equiv f_2(f_1 (\bm{x}^t)), \ \bm{x}^t \in \bm{B}^N \equiv L_D
\end{equation}
Fig. \ref{fg4} shows the Cmap for the 6-dimensional PBNN in Example 2: 
composing the 1st map in Fig. \ref{fg4}(a) and 2nd map in Fig. \ref{fg4}(b), 
we obtain the Cmap in Fig. \ref{fg4} (c). 
The 1st map visualizes dynamics of the 6-dimsnsional SBNN in Fig. \ref{fg2} (a),  
the 2nd map visualizes global permutation connection P231465, 
and the Cmap visualizes dynamics of the 6-dimensional PBNN in Fig. \ref{fg2} (b). 
The Cmap exhibits multiple BPOs including the BPO with period 12 as shown in Fig. \ref{fg4} (d). 
This figure suggests that the global permutation connection is effective to realize a variety of BPOs. 
The BPO, BPP, and EPP of the Cmap are defined by replacing the 1st map $f_1$ with the Cmap $f$ in Definition 1. 

\subsection{Tool 2: Feature quantities}
In order to classify PBNNs by their BPOs, we present the Tool 2: two simple feature quantities. 
As shown in the Cmap in Fig. \ref{fg4} (d), the PBNN can have multiple BPOs and 
exhibits one of them depending on the initial points. 
Since analysis of multiple BPOs is not easy, we consider one BPO with the maximum period (MBPO). 
In Example 2, the BPO with period 12 is the MBPO as shown in Fig. \ref{fg4}(d). 
The first feature quantity is defined by 
\begin{equation}
\alpha = \mbox{(The period of MBPO)}/2^N, \  \frac{1}{2^N} \le \alpha \le 1. 
\end{equation}
As $\alpha$ increase, the MBPO approaches the M-sequence. 
Roughly speaking, the quantity $\alpha$ evaluates complexity of the MBPO. 

The second feature quantity is defined by
\begin{equation}
%\begin{array}{c}
\beta = \mbox{(The number of initial points falling into the MBPO)}/2^N, \ 
\alpha \le \beta \le 1.
%\end{array}
\end{equation}
Note that the initial points consist of BPPs in the MBPO and EPPs to the MBPO: $\alpha \le \beta$. 
The difference $(\beta - \alpha)$ is proportional to the number of EPPs that evaluates stability of the MBPO. 
In order not to disturb existence of EPPs, the period should not be too long, e.g., $\alpha \le 1/3$. 
Note that the stability corresponds to bit error correction capability in digital dynamical systems and 
corresponds to basin of attraction in analog dynamical systems \cite{ott}.  
As $\beta$ increases, the stability becomes stronger. 
If multiple MBPOs exist, one MBPO with larger $\beta$ is adopted.

\begin{figure}[tb]
\centering
\includegraphics[width=0.5\columnwidth]{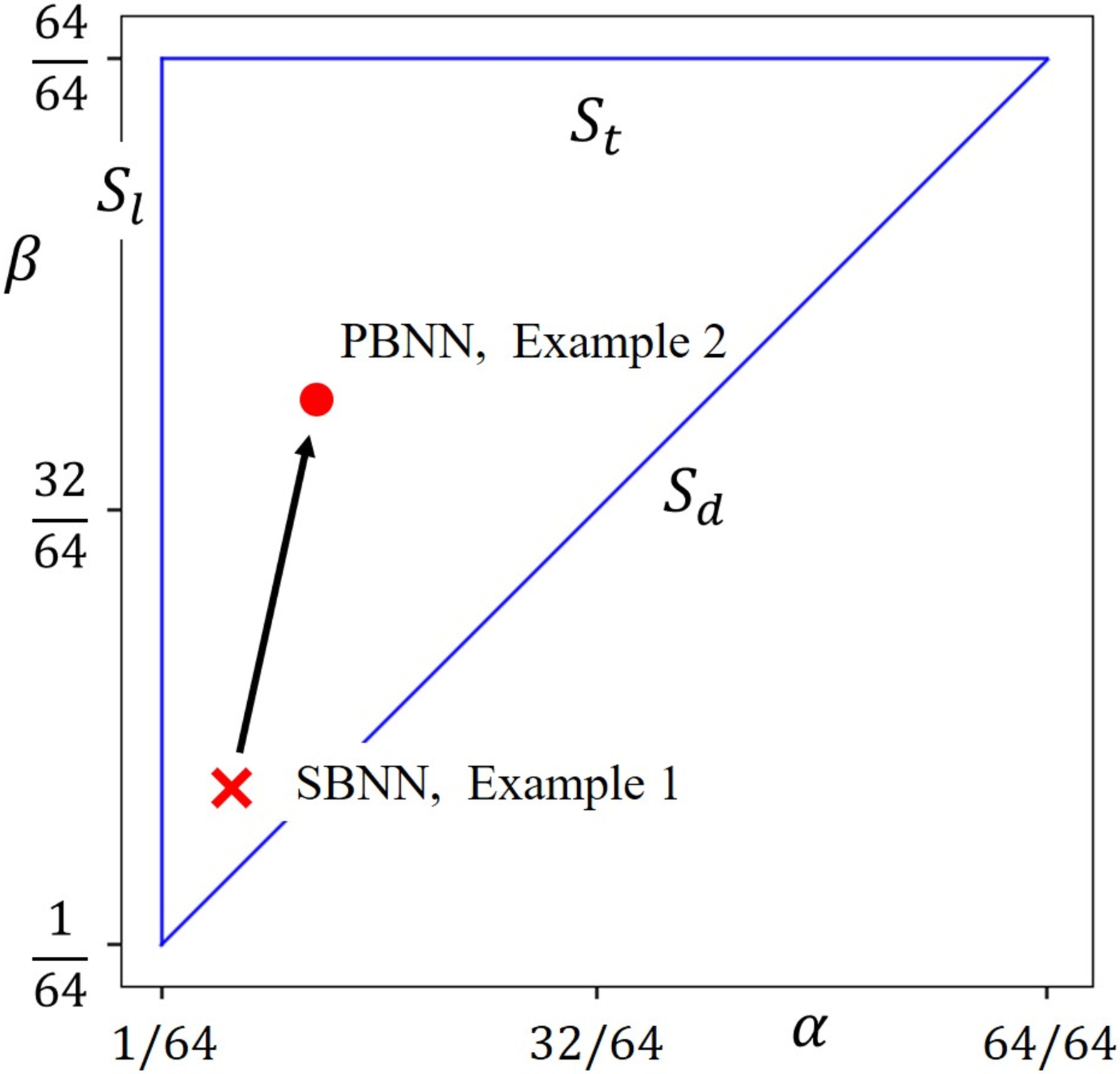}
\caption{Feature plane and three segments $S_d$, $S_t$, and $S_l$. 
Red cross: SBNN, Example 1 ($\alpha=6/64$, $\beta=12/64$). 
Red circle: PBNN, Example 2 ($\alpha=12/64$, $\beta=40/64$). 
}
\label{fg5} 
\end{figure}

Using the two feature quantities, we construct a feature plane as shown in Fig. \ref{fg5}. 
In the figure, 
SBNN in Example 1 is represented by a point $(\alpha, \beta)=(6/64, 12/64)$ and 
PBNN in Example 2 is represented by a point $(\alpha, \beta)=(12/64, 40/64)$. 
The permutation (P231465) increases both $\alpha$ and $\beta$. 

As stated earlier, as $\alpha$ increases, period of the MBPO becomes longer (more complex). 
As $(\beta-\alpha)$ increases, stability of the MBPO becomes stronger. 
A PBNN is represented by its MBPO and the MBPO gives one point in the feature plane. 
The point exists in a triangle surrounded by the three segments each of which has the following meaning. 
\begin{itemize}
\item $S_d (\alpha=\beta)$: There exists no EPP to the MBPO and no transient phenomenon to the MBPO. 
\item $S_t (\beta=1)$: All the initial points fall into the MBPO and the PBNN has unique BPO=MBPO. 
\item$S_l (\alpha=1/2^N)$: The MBPO is a fixed point (BPO with period 1). 
\end{itemize}
\noindent
The feature plane is useful for basic classification of PBNNs.

\section{Numerical analysis}
\label{nanalysis}

Using the Tool 1 and Tool 2, we analyze PBNNs. 
Since analysis of $N$-dimensional PBNNs is not easy, we consider 6-dimensional PBNNs in this paper. 
The case $N=6$ is convenient in the following three points. 
First, the case $N=6$ is applicable to engineering systems such as 
hexapod walking robots \cite{cnna} and power converters with 6 switches \cite{dbnn1}. 
Second, in the 6-dimensional PBNNs, the number of points in the parameter space (parameter values) is $8 \times 6! = 5760$ and the brute force attack is possible to calculate Cmaps and feature quantities. 
The 5760 points is much smaller than 
$3^{6^2}$ points in the parameter space of the 6-dimensional DBNNs in Eq. (\ref{dbnn}). 
An effective calculation method for BPOs and EPPs can be found in \cite{nolta}.  
As $N$ increases, the number of elements in the domain $\bm{B}^N \equiv L_N$ 
increases exponentially and we must use an evolutionary computation algorithm etc. 
Third, the Cmap consists of $2^N=2^6$ points and is visible.  
As $N$ increases, the number of points increases exponentially and the Cmap is hard to inspect. 

\begin{figure}[htb]%
\centering
\includegraphics[width=0.85\columnwidth]{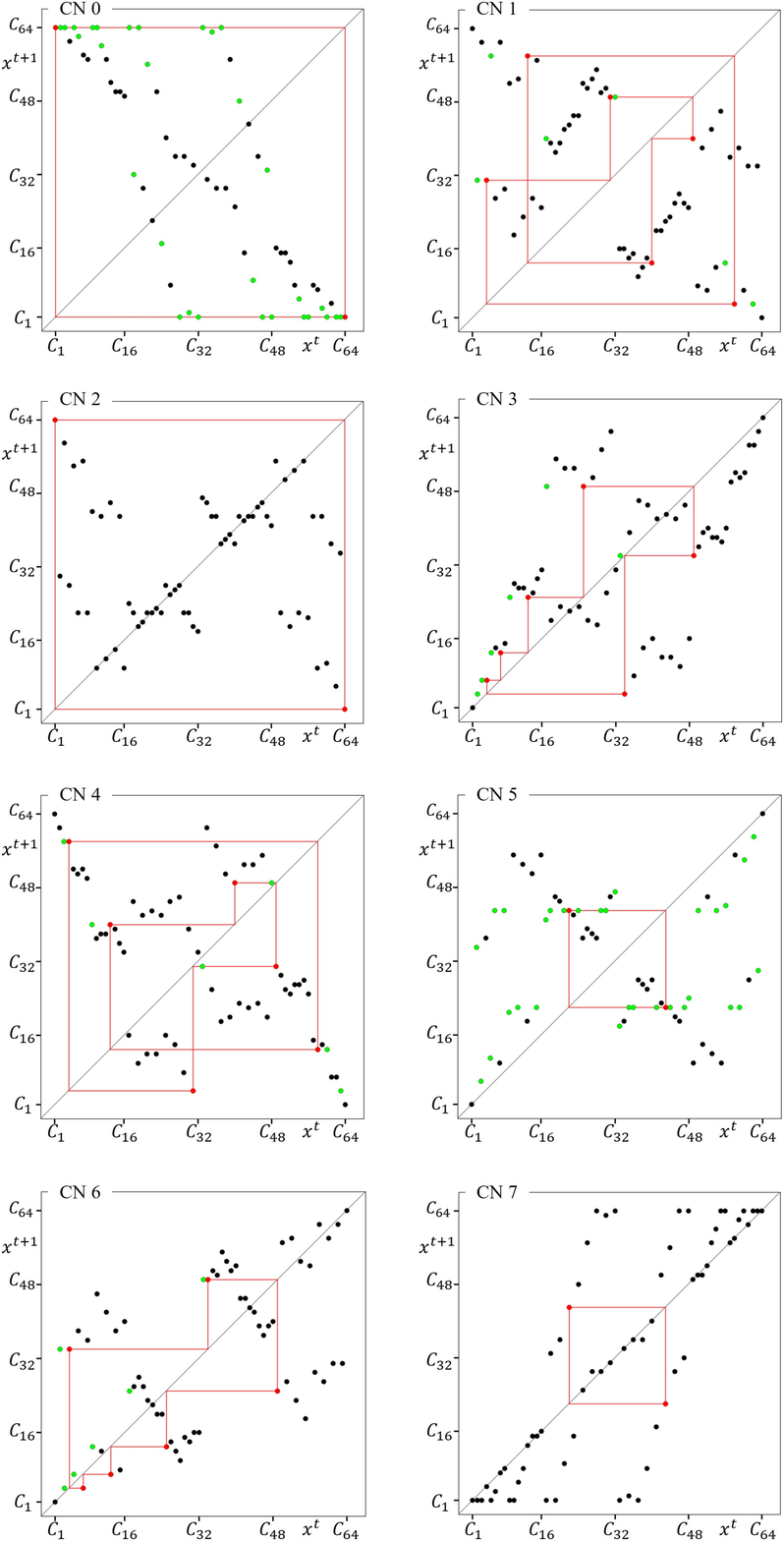}
\caption{1st map of 8 SBNNs 
with MBPO (red) and EPPs (green). 
%CN0: MBPO with period 2. 
%CN1: period  6. 
%CN2: period  2. 
%CN3: period  6. 
%CN4: period  6. 
%CN5: period  2. 
%CN6: period  6. 
%CN7: period  2.
} 
\label{fg6} 
\end{figure}

Fig. \ref{fg6} shows 1st maps with MBPOs for all the 8 SBNNs. 
As stated earlier, the SBNNs correspond to 8 ECAs governed by 8 rules of $\lambda$-parameter 0.5. 
In CN6, the MBPO with period 6 is applicable to control signal of hexapod robots \cite{cnna}. 
Using the two feature quantities, MBPOs of the 8 SBNNs are evaluated: 
%
%\begin{equation}
\[
\begin{array}{lll}
\mbox{CN0}: (w_a, w_b, w_c) = (-1, -1, -1)  & \alpha = 2/64  & \beta = 32/64\\ 
\mbox{CN1}: (w_a, w_b, w_c) = (-1, -1, +1)  & \alpha = 6/64  & \beta = 12/64\\ 
\mbox{CN2}: (w_a, w_b, w_c) = (-1, +1, -1)  & \alpha = 2/64  & \beta = 2/64\\ 
\mbox{CN3}: (w_a, w_b, w_c) = (-1, +1, +1)  & \alpha = 6/64  & \beta = 12/64\\ 
\mbox{CN4}: (w_a, w_b, w_c) = (+1, -1, -1)  & \alpha = 6/64  & \beta = 12/64\\ 
\mbox{CN5}: (w_a, w_b, w_c) = (+1, -1, +1)  & \alpha = 2/64  & \beta = 32/64\\ 
\mbox{CN6}: (w_a, w_b, w_c) = (+1, +1, -1)  & \alpha = 6/64  & \beta = 12/64\\ 
\mbox{CN7}: (w_a, w_b, w_c) = (+1, +1, +1)  & \alpha = 2/64  & \beta = 2/64
\end{array}
\]
%\end{equation}
%
Dynamics of all the 8 SBNNs are visualized in the 8 1st maps in Fig. \ref{fg6} 
corresponding to all the 8 points in the parameter space. 
The SBNNs exhibit MBPOs with period 6 and MBPOs with period 2. 

\begin{figure}[htb]
\centering
\includegraphics[width=0.85\columnwidth]{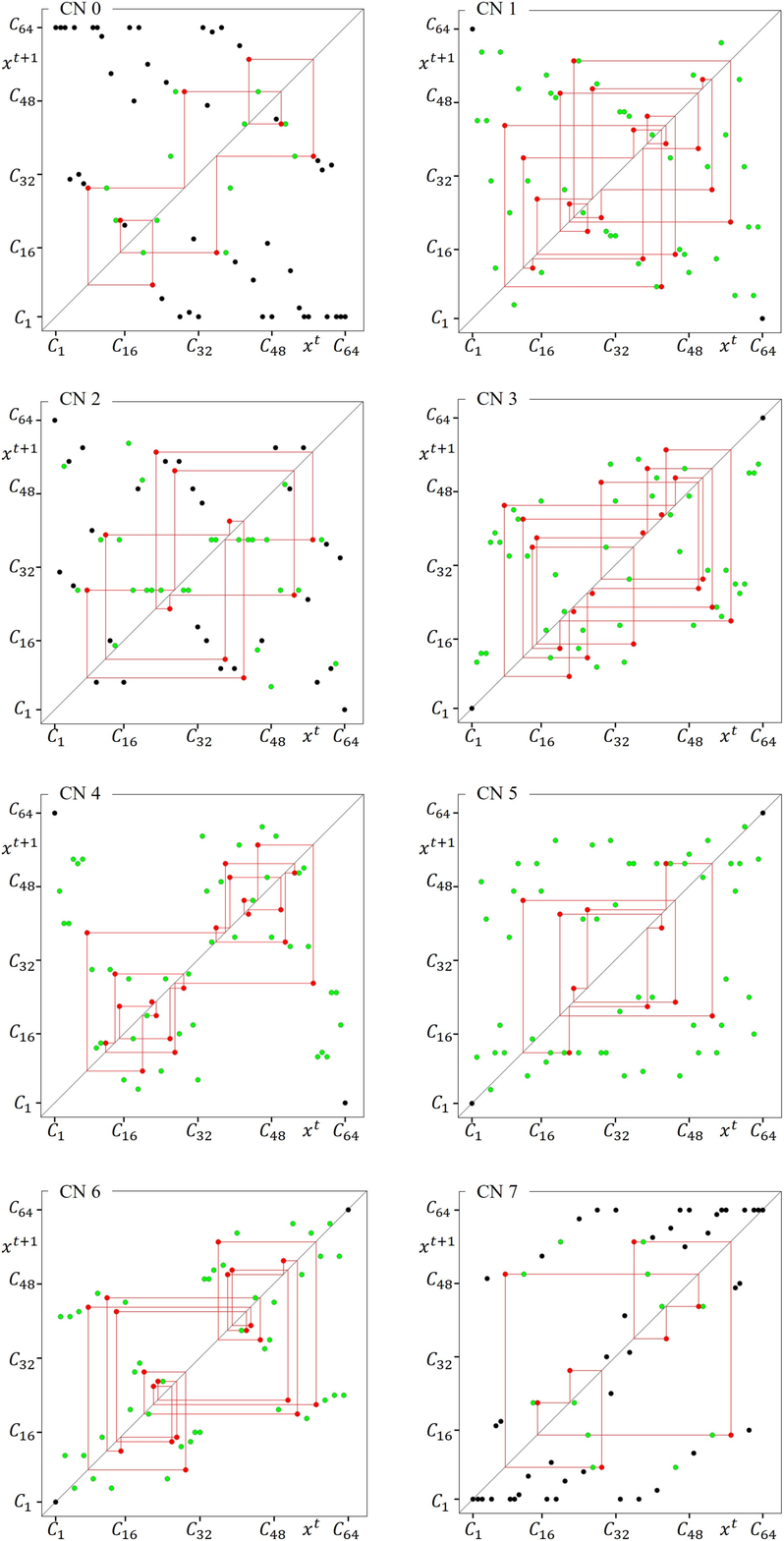}
\caption{Cmap of PBNNs with MBPO (red) and EPPs (green). 
CN0: P513246. 
CN1: P413625. 
CN2: P524361. 
CN3: P315462. 
CN4: P254136. 
CN5: P461253. 
CN6: P126354. 
CN7: P651324. 
}
\label{fg7} 
\end{figure}

Applying global permutation connection to each SBNN of local binary connection, 
we obtain PBNNs that can generate MBPO with longer period than that of SBNN. 
Fig. \ref{fg7} shows Cmaps of a typical PBNN for each SBNN. 
Using the two feature quantities, MBPOs of the 8 PBNNs are evaluated:  
\begin{equation}
\begin{array}{llll}
\mbox{CN0}: \mbox{P513246}  & \alpha = 8/64   & \beta = 20/64\\ 
\mbox{CN1}: \mbox{P413625}  & \alpha = 20/64  & \beta = 62/64\\ 
\mbox{CN2}: \mbox{P524361}  & \alpha = 10/64  & \beta = 36/64\\ 
\mbox{CN3}: \mbox{P315462}  & \alpha = 20/64  & \beta = 62/64\\ 
\mbox{CN4}: \mbox{P254136}  & \alpha = 20/64  & \beta = 62/64\\ 
\mbox{CN5}: \mbox{P461253}  & \alpha = 10/64  & \beta = 62/64\\ 
\mbox{CN6}: \mbox{P126354}  & \alpha = 20/64  & \beta = 62/64\\ 
\mbox{CN7}: \mbox{P651324}  & \alpha = 8/64   & \beta = 20/64\\ 
\end{array}
\label{fq}
\end{equation}
Dynamics of the 8 PBNNs are visualized in 8 Cmaps in Fig. \ref{fg6}.  
In these examples, four PBNNs give the MBPO with the maximum period 20 ($\alpha=20/64$) and 
the strongest stability ($\beta=62/64$). 
In each PBNN, the global permutation connection realizes longer period than each SBNN.

\begin{figure}[htb]
\centering
\includegraphics[width=0.85\columnwidth]{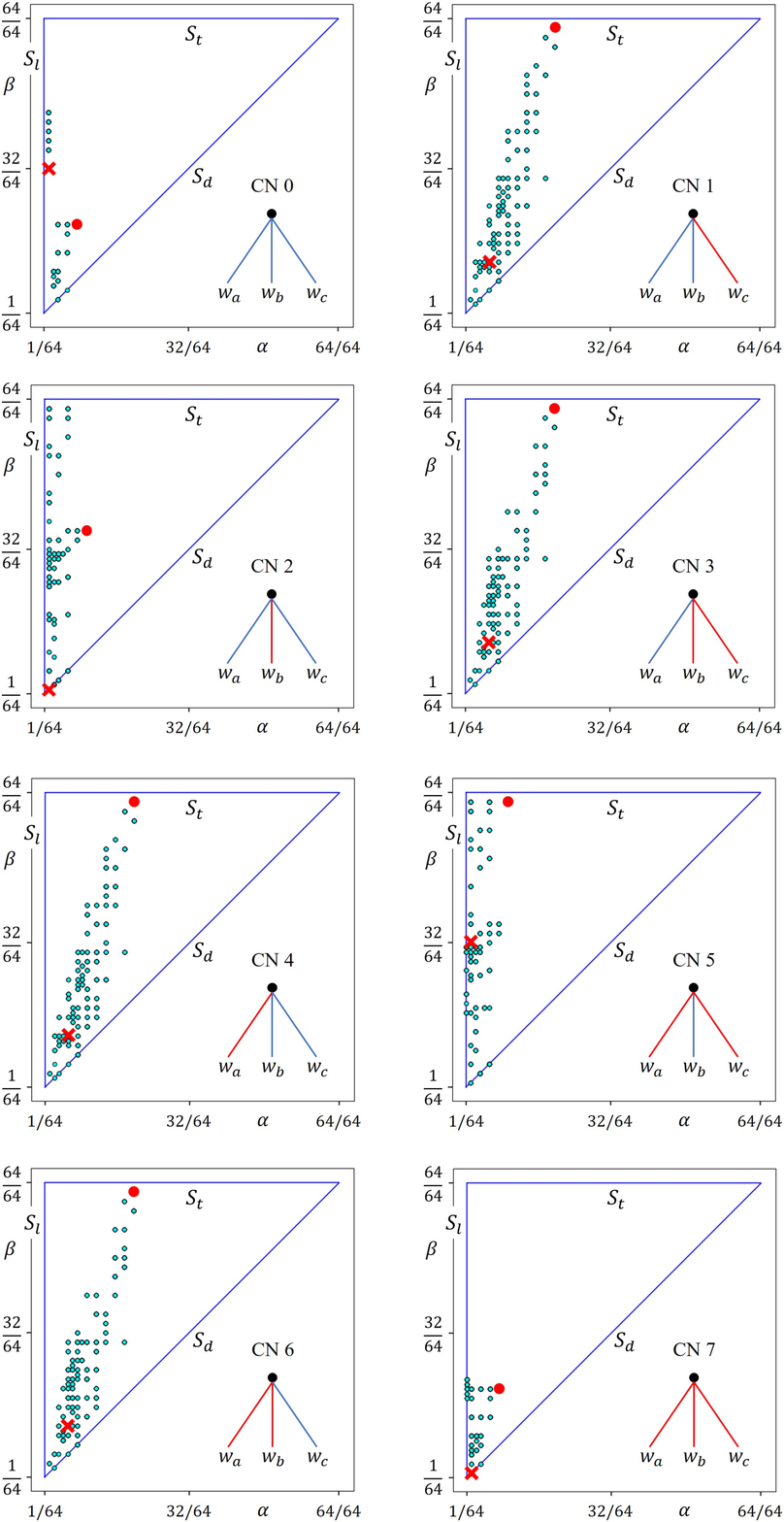}
\caption{Feature planes.  
Red cross: SBNN, Red circle: PBNN in Fig. \ref{fg7}. 
CN0(19 points). CN1(79). CN2(49). CN3(78). 
CN4(79).        CN5(53). CN6(78). CN7(26). 
}
\label{fg8} 
\end{figure}

Next, applying the brute force attack to all $6!=720$ PBNNs for each SBNN, 
we have calculated the two feature quantities. 
The results are summarized in the feature plane in Fig. \ref{fg8}. 
In each feature plane, the red cross corresponds to the SBNN of the 1st map in Fig. \ref{fg6} and  
the red circle corresponds to the PBNN of the Cmap in Fig. \ref{fg7}. 
The PBNN of the red circle has the MBPO with the longest period ($\alpha$) and 
the strongest stability ($\beta$) for the longest period. 
We can see that the global permutation connection is effective to increase the period and stability of MBPOs. 
Note that points of CN1 and CN4 (respectively, CN3 and CN6) are the same and the points are distributed in wide range: the global permutation connection is effective to generate a variety of MBPOs. 
In CN2 and CN5, $\alpha$ distributes in a narrow range whereas $\beta$ distributes in a wide range. 
In CN0 and CN7, both $\alpha$ and $\beta$ distributes in a relatively narrow range: 
the local binary connection $w=(-1,-1,-1)$ or $(+1,+1,+1)$ seems not to make a variety of BPOs. 
These results from 6-dimansional PBNNs suggest huge variety of BPOs and 
hardness of analysis in the higher-dimensional PBNNs. 
Note that, in Fig. \ref{fg8}, the parameter space of the PBNNs is classified into 8 subspaces (8 CNs). 
On each parameter subspace, the two feature quantities are distributed and the typical values are red marked. 
These results provide basic information to consider relationship between 
the feature quantities and the parameters. 
The parameters are represented by the permutation identifier and connection number (CN).

\clearpage
It should be noted that, on each feature plane, the number of points is less than the number of PBNNs. 
For example, on the feature plane of CN6, 78 points exist for 720 PBNNs. 
Multiple PBNNs give the same value of the two future quantities $\alpha$ and $\beta$. 
For more detailed classification, novel feature quantities are necessary.

\section{Tool 3: FPGA based hardware}
In order to realize engineering applications, we present the Tool 3: an FPGA based hardware prototype. 
The FPGA is an integrated circuit designed to be configured by a designer after manufacturing.   
The advantages include high speed/precision operation and high degree of integration \cite{fpga} \cite{dsn}. 
We have designed the hardware and have measured typical BPOs in the following environment. 
%\noindent $\bullet$ 
%
\begin{itemize}
\item Verilog version: Vivado 2018.2 platform (Xilinx).
\item FPGA: BASYS3 (Xilinx Artix-7 XC7A35T-ICPG236C).
\item Clock: 10[MHz]. The default frequency 100[MHz] is divided for clear measurements. 
\item Measuring instrument: ANALOG DISCOVERY2. 
\item Multi-instrument software: Waveforms 2015
\end{itemize}
Algorithms \ref{alg1} and \ref{alg2} show SystemVerilog codes for hardware design. 
In the Algorithm \ref{alg1}, as a connection number CN (equivalent to a rule number RN) is given, 
the corresponding Boolean function is obtained and the SBNN is realized.   
This design strategy is based on the Boolean functions from three input to one output 
and is applicable to implementation of ECAs for all the 256 rules. 
In Algorithm \ref{alg2}, as a permutation identifier is given, the global permutation connection is realized. 
Applying these two algorithms, a desired PBNN is realized. 
Varying the CN and permutation identifier, we obtain a desired PBNN in $8 \times N!$ PBNNs. 

\begin{figure}[tb]
\centering
\includegraphics[width=0.3\columnwidth]{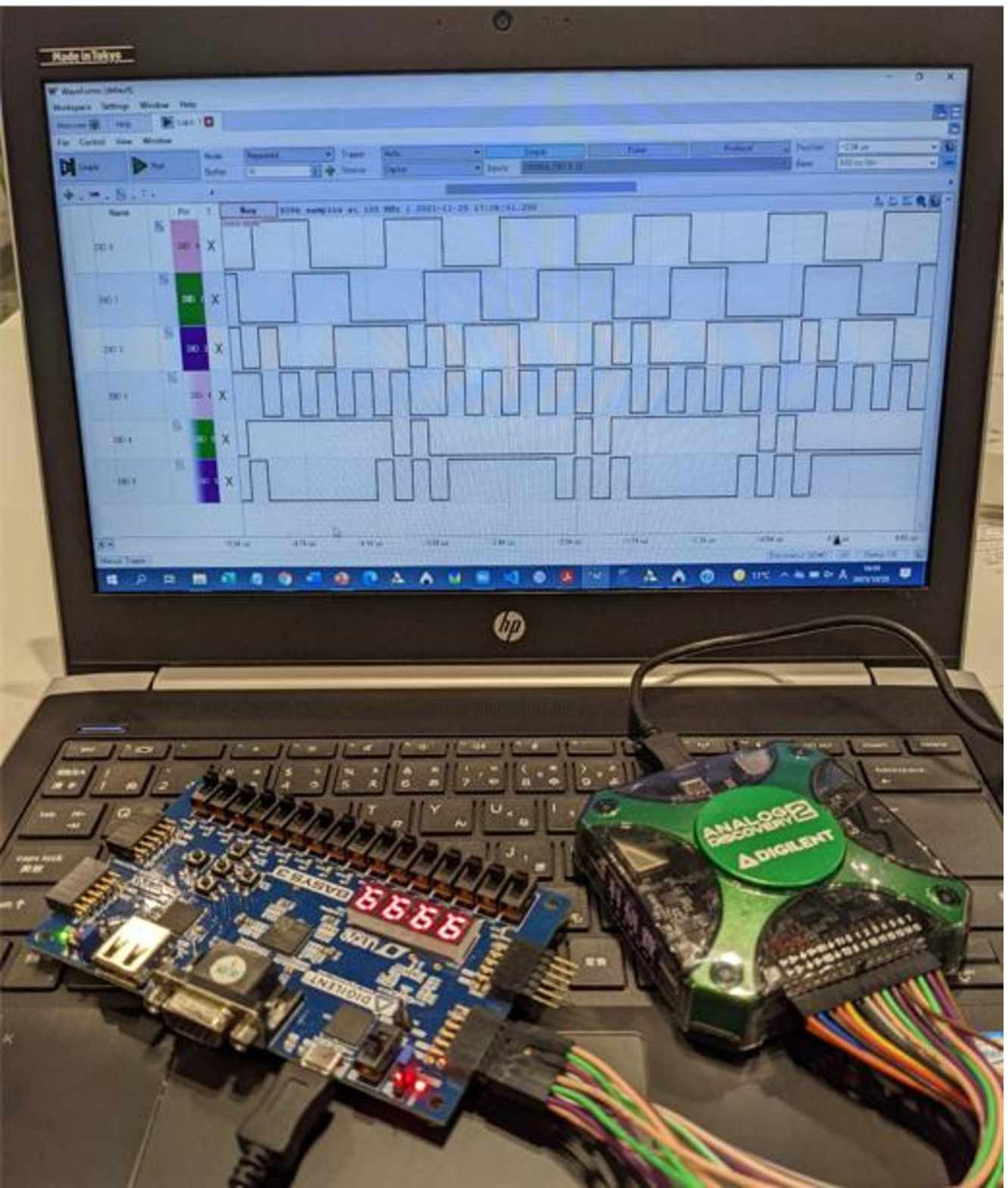}
\caption{
Experimental setup: FPGA and Analog discovery.
}
\label{fg9}
%\end{figure}

%\begin{figure}[htb]
\centering
\includegraphics[width=0.8\columnwidth]{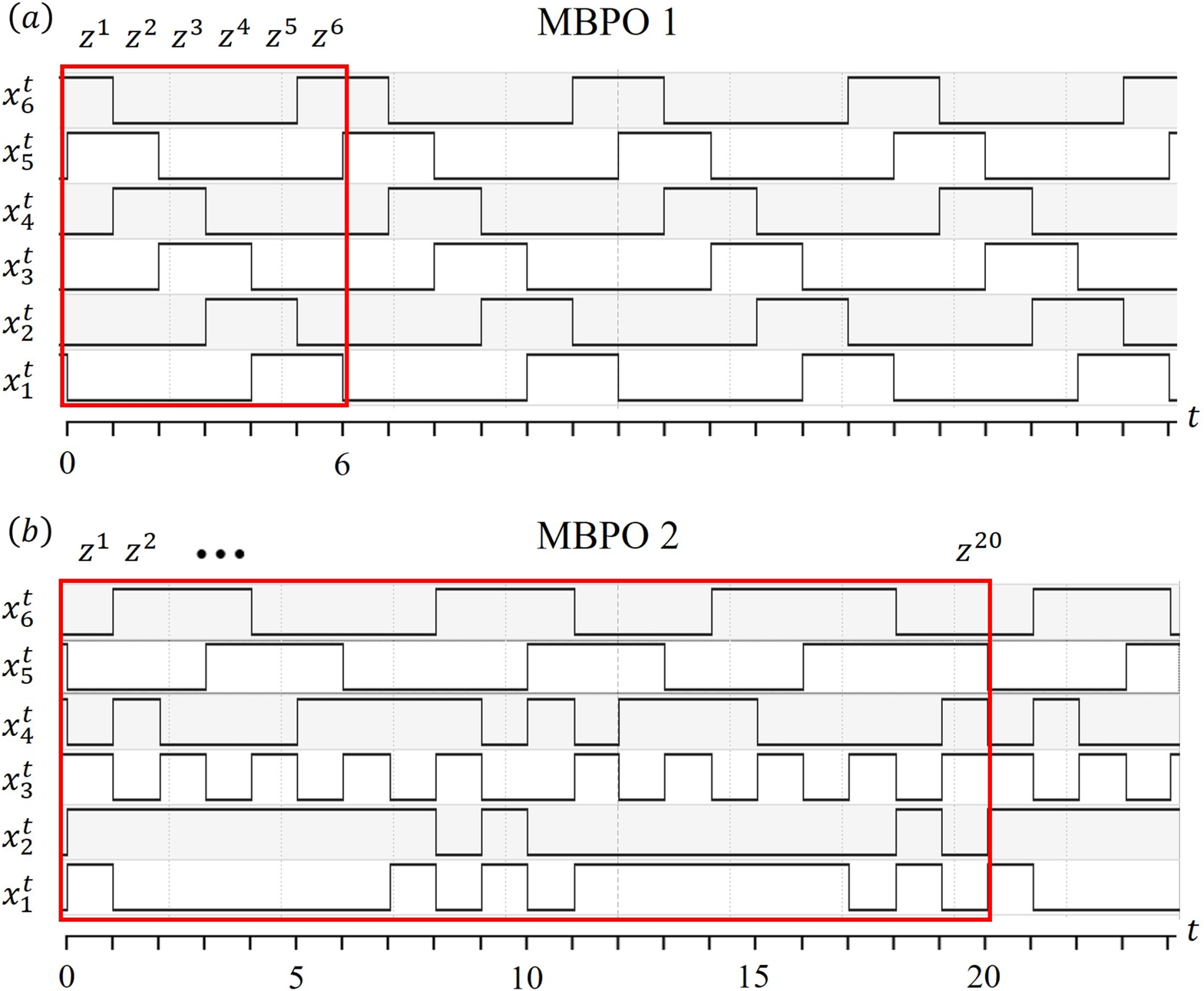}
\caption{
Measured waveforms of MBPOs in the FPGA board. 
(a) MBPO1 with period 6 from SBNN CN6. 
(b) MBPO2 with period 20 from PBNN CN6 P126354. 
}
\label{fg10}
\end{figure}

Using the algorithms, the SBNN and PBNN are implemented on the FPGA board. 
Fig. \ref{fg10}(a) shows measured waveform of MBPO1 with period 6 form the SBNN in Fig. \ref{fg2} (a). 
The MBPO1 is applicable to control signal of hexapod walking robots 
and the control strategy can be found in \cite{cnna}. 
%Although this paper does not aim at application to hexapod robots, 
Fig. \ref{fg10}(b) shows measured waveform of MBPO2 with period 20 corresponding to the Cmap in Fig. \ref{fg7} CN6. 
The MBPO2 has the longest period that is impossible in 6-dimensional SBNNs. 
This FPGA prototype provides basic information to realize engineering applications of 6 or higher dimensional PBNNs.

\begin{algorithm}
	\caption{SBNN}
	\label{alg1}
	\begin{algorithmic}
		\STATE module SBNN(parameter $N=6$)
		\STATE (output reg $x[1:N]$,
		\STATE reg $[1:N] x^{t+1}$);
                \STATE wire $rule0,rule1, \cdots, rule7;$
                %\STATE integer $w[0:7]=[0, 0, 1, 0, 1, 0, 1, 1];$
		\STATE parameter RN = 8'b212 ; 
		\hfill // Rule Number 
		\STATE genvar $j$;
		\STATE for$(j=1; j<=N; j=j+1)$begin
		\STATE \ \ \ \ assign $ rule0=RN[0]^{\ast}($\textasciitilde$x[j-1] \ \&\ $\textasciitilde$x[j] \ \& \ $\textasciitilde$x[j+1]);$
		\hfill // Boolean function
		\STATE \ \ \ \ assign $ rule1=RN[1]^{\ast}($\textasciitilde$x[j-1] \ \&\ $\textasciitilde$x[j] \ \& \ x[j+1]);$
		\STATE \ \ \ \ assign $ rule2=RN[2]^{\ast}($\textasciitilde$x[j-1] \ \&\ x[j] \ \& \ $\textasciitilde$x[j+1]);$
		\STATE \ \ \ \ assign $ rule3=RN[3]^{\ast}($\textasciitilde$x[j-1] \ \&\ x[j] \ \& \ x[j+1]);$
		\STATE \ \ \ \ assign $ rule4=RN[4]^{\ast}(x[j-1] \ \&\ $\textasciitilde$x[j] \ \& \ $\textasciitilde$x[j+1]);$
		\STATE \ \ \ \ assign $ rule5=RN[5]^{\ast}(x[j-1] \ \&\ $\textasciitilde$x[j] \ \& \ x[j+1]);$
		\STATE \ \ \ \ assign $ rule6=RN[6]^{\ast}(x[j-1] \ \&\ x[j] \ \& \ $\textasciitilde$x[j+1]);$
		\STATE \ \ \ \ assign $ rule7=RN[7]^{\ast}(x[j-1] \ \&\ x[j] \ \& \ x[j+1]);$
                \STATE \ \ \ \ assign $x^{t+1}[j] = (rule0)|(rule1)|~...~|(rule7);$%assign $x^{t+1}[j] = (rule0)|(rule1)|(rule2)|(rule3)|(rule4)|(rule5)|(rule6)|(rule7);$
		\STATE end
		\STATE endmodule
	\end{algorithmic}
\end{algorithm}

\begin{algorithm}
	\caption{PBNN}
	\label{alg2}
	\begin{algorithmic}
		\STATE module PBNN(parameter $N=6$)
		\STATE (input clk$,$load$,$rst$,$
                \STATE input $i[1:N],$
		\STATE output reg $x[1:N]$);
		\STATE reg $[1:N] x^{t+1}$;
		\STATE integer $k$;
		\STATE integer $y[1:N]=[1, 2, 6, 3, 5, 4];$ \hfill // Permutation identifier P126354
		\STATE always $@$(posedge clk)begin
		\STATE \ \ if(load$==1$)begin
		\STATE \ \ \ \ for$(k=1; k<=N; k=k+1)$begin
		\STATE \ \ \ \ \ \ $x[k] = i[k];$ \hfill // Initial condition
		\STATE \ \ \ \ end
		\STATE \ \ end else if(rst$==1$)begin
		\STATE \ \ \ \ for$(k=1; k<=N; k=k+1)$begin
		\STATE \ \ \ \ \ \ $x[k] = 0;$
		\STATE \ \ \ \ end
		\STATE \ \ end else begin
		\STATE \ \ \ \ for$(k=1; k<=N; k=k+1)$begin
		\STATE \ \ \ \ \ \ $x[k] = x^{t+1}[y[k]];$ \hfill // Permutation
		\STATE \ \ \ \ end
		\STATE \ \ end
		\STATE end
		\STATE SBNN SBNN$(x,x^{t+1})$;
		\STATE endmodule
	\end{algorithmic}
\end{algorithm}

\section{Conclusions}
\label{concl}
The PBNNs and the three tools are presented in this paper: 
the Cmap for visualization of the dynamics,  
the two feature quantities for basic classification of PBNNs, and 
the FPGA based hardware prototype for engineering applications. 
Using the tools, all the 6-dimensionl PBNNs have been analyzed and 
typical MBPOs are confirmed experimentally. 
In our future works, we should consider various problems including the following. 
%In the two feature quantities $\alpha$ and $\beta$, multiple PBNNs give the same value. 
1. Relation between permutation identifier and feature quantities. 
2. Application of the global permutation to ECAs of 256 rules and analysis of the dynamics. 
3. Novel feature quantities for more detailed classification of PBNNs. 
4. Evolutionary algorithms for analysis and synthesis of higher-dimensional PBNNs. 
5. Engineering applications of the FPGA based hardware. 
%6-10D for hexapod robots, Higher D for signal processing. 
%Simple topology with rich BPO: picking/growing topology is important/hard \cite{top}.  

\clearpage

\medskip
% The data information below will be filled by AIMS editorial staff
%Received xxxx; revised xxxx.
\medskip

\end{document}